\documentclass[10pt,a4paper]{article}
\usepackage[T1]{fontenc}
\usepackage[utf8]{inputenc}
\usepackage{lmodern}
\usepackage[margin=25mm]{geometry}
\usepackage{amsmath,amssymb,amsthm}
\usepackage{mathtools}
\usepackage{enumitem}
\usepackage{microtype}
\usepackage{hyperref}
\hypersetup{hidelinks,pdfauthor={Shuyuan Liang}}
\theoremstyle{plain}
\newtheorem{theorem}{Theorem}[section]
\newtheorem{lemma}[theorem]{Lemma}
\newtheorem{proposition}[theorem]{Proposition}

\theoremstyle{remark}
\newtheorem{remark}{Remark}
\numberwithin{equation}{section}

\newcommand{\Pp}{\mathbf P}
\newcommand{\E}{\mathbf E}
\newcommand{\R}{\mathbb R}
\newcommand{\ind}{\mathbf 1}
\newcommand{\eqdist}{\mathrel{\stackrel{\mathrm d}{=}}}
\newcommand{\dto}{\xrightarrow{\mathrm d}}
\newcommand{\Phibar}{\overline\Phi}
\DeclareMathOperator{\Var}{Var}
\DeclareMathOperator{\Law}{\mathcal L}

\newenvironment{manuscriptproof}[1][Proof]
 {\par\smallskip\noindent\textbf{#1.}\ \ignorespaces}
 {\unskip\nobreak\hfill\qed\par\smallskip}
\allowdisplaybreaks[1]

\title{Normal approximation and self-normalized moderate deviations
for the minimal random walk with random step sizes}
\author{Shuyuan Liang\thanks{Corresponding author. Email:
\href{mailto:liangshuyuan0501@foxmail.com}{\nolinkurl{liangshuyuan0501@foxmail.com}}}\\[0.5ex]
\small School of Mathematics and Statistics, Northeastern University at Qinhuangdao,\\
\small Qinhuangdao, {\rm066004}, Hebei, China.}
\date{}

\begin{document}
\maketitle

\begin{abstract}
We study the minimal random walk with independent positive random step
sizes. By separating the movement count from the centered step-size
noise, we obtain normal approximations and self-normalized Cram\'er-type
moderate deviations under finite-moment assumptions. For the noise,
conditioning on the movement count yields Kolmogorov, Wasserstein and
Zolotarev bounds. A polynomial moderate-deviation range holds when a
recalled stop can lead to movement. At the boundary where this
probability is zero, we obtain a logarithmic range and construct a
bounded nonsymmetric example excluding a general polynomial range.
For the full position, a weighted martingale representation gives
limit laws in the different regimes, nonuniform normal approximations
and moment-distance bounds without logarithmic losses. Variance
concentration also yields self-normalized moderate deviations in the
diffusive and critical regimes.
\end{abstract}

\begingroup\small

\par\endgroup

\section{Introduction}\label{sec:introduction}
\quad
The minimal random walk (MRW), introduced by Harbola, Kumar and
Lindenberg \cite{Harbola2014}, is a random walk with memory whose
increments take the values zero and one. At time $n+1$, an index is
chosen uniformly from $\{1,\ldots,n\}$. If the recalled increment is
one, the next increment is one with probability $p$; if it is zero,
the next increment is one with probability $q$. Otherwise the walk
remains at its current position. Writing $I_k$ for the increments,
with $I_1=1$, the movement count satisfies
\[
 K_n=\sum_{k=1}^n I_k,\qquad
 \Pp(I_{n+1}=1\mid I_1,\ldots,I_n)=q+(p-q)K_n/n.
\]
For $q>0$ and $p<1$, the walk is called diffusive, critical and
superdiffusive according to $p-q<1/2$, $p-q=1/2$ and $p-q>1/2$,
respectively.

The asymptotic behavior of the MRW has been studied by several
authors. Coletti, Gava and de Lima \cite{Coletti2019} established
laws of large numbers and central limit theorems. Bercu and
V\'azquez Guevara \cite{Bercu2022} developed a martingale approach
to the three regimes. Their results include almost sure limit
theorems and Gaussian fluctuations around the random limit in the
superdiffusive regime. Miyazaki and Takei \cite{Miyazaki2020}
studied the boundary case $q=0$, where recalling a zero increment
forces another stop. For $0<p<1$, the movement count then has
order $n^p$ and its normalization converges to a Mittag--Leffler
random variable.

An extension of this model attaches independent positive random
step sizes to the movements. Let $(Z_k)_{k\ge1}$ be i.i.d.,
independent of the movement process, with mean one. The position
is $H_n=\sum_{k=1}^n I_kZ_k$, and its centered step-size noise is
$R_n=H_n-K_n$. Thus
\[
 H_n-\E H_n=(K_n-\E K_n)+R_n.
\]
The normalization of the full position must account for both terms.
Conditional on the movement process, $R_n$ is a sum of $K_n$ independent centered
step sizes. The full position also contains the fluctuations of
$K_n$, and consequently retains the phase transition of the MRW.

In this paper, we study normal approximation and self-normalized
Cram\'er-type moderate deviations for both quantities. For $R_n$,
we combine the independent-sum estimates of
\cite{Petrov1975,Shao2005,Jing2003,Rio2009} with negative moments
and lower-tail estimates for $K_n$; see also \cite{deLaPena2009}
for the general theory of self-normalization. This gives Kolmogorov
bounds under finite moments and Wasserstein and Zolotarev bounds
for bounded step sizes. For the self-normalized
statistic, centering and variance matching are included when
required by the order of the Zolotarev distance. Under a nonzero
third centered moment, a matching lower bound shows that the
third-order Zolotarev rate for the variance-standardized noise is
attained. No symmetry of the centered step sizes is assumed.

The moderate-deviation range depends on the movement count.
When $q>0$, the usual polynomial range for independent sums is
preserved. When $q=0$, small movement counts yield a logarithmic
sufficient range, for which we give an explicit endpoint with
leading term $\sqrt{2p\log n}$. A bounded nonsymmetric example
shows that a polynomial range cannot hold for all step-size laws
under these assumptions.

For $H_n-\E H_n$, we use a weighted martingale array containing
both the movement innovations and the centered step-size terms.
The movement limit theorems yield the corresponding limits for
the full position, including a Gaussian fluctuation about a
random center in the superdiffusive regime. In the diffusive and
critical regimes, variance concentration and moment estimates
allow us to apply the martingale bounds of
\cite{HaeuslerJoos1988,FanShao2018,FanShao2025}.
Integration of the nonuniform normal approximation gives
Wasserstein and Zolotarev bounds under finite moments, without
additional logarithmic factors. At $q=0$, the variance-standardized
full position instead converges to a non-Gaussian limit.

Section~\ref{sec:main} states the model and main results.
Section~\ref{sec:clock} collects the preliminary estimates, and
Sections~\ref{sec:normal-proofs} and \ref{sec:md-proofs} contain
the proofs for the centered step-size noise.
Section~\ref{sec:full-proofs} treats the full position, and
Section~\ref{sec:conclusion} concludes the paper.

Throughout the paper, $c$ and $C$, possibly with indices, denote
positive constants depending only on fixed model parameters, the
step-size distribution and moment indices. Their values may vary
from line to line. For positive sequences, $u_n\asymp v_n$ means
that their ratio is bounded above and away from zero for all
sufficiently large $n$, and $u_n\sim v_n$ means that their ratio
tends to one. Let $G$ have the standard normal distribution,
with distribution function $\Phi$, and put $\Phibar=1-\Phi$.
We write $U\eqdist V$ for equality in distribution and
$\|U\|_s=(\E|U|^s)^{1/s}$ for $s\ge1$.
Uniformity for $0\le x=o(h_n)$ means uniformity on $[0,u_n]$
for every nonnegative deterministic sequence $u_n=o(h_n)$.

\section{The model and main results}\label{sec:main}

\subsection{The movement process and random step sizes}

\quad
The MRW with random step sizes is defined as follows. Let $I_1=1$ and
\begin{equation}\label{eq:mrw}
 K_n=\sum_{k=1}^n I_k,\qquad \mathcal F_n=\sigma(I_1,\ldots,I_n).
\end{equation}
The minimal random walk is determined by
\begin{equation}\label{eq:transition}
 \Pp(I_{n+1}=1\mid\mathcal F_n)=q+(p-q)K_n/n,
 \qquad p,q\in[0,1].
\end{equation}

We consider the following two parameter sets:
\begin{equation}\label{eq:regimes}
 \begin{split}
 \text{(A)}\quad &q=0,\qquad p\in(0,1),\\
 \text{(B)}\quad &q\in(0,1],\qquad p\in[0,1].
 \end{split}
\end{equation}
The deterministic initialization avoids the absorbing all-stop path
that is possible in case~(A) with a random initial indicator.

Let $(Z_k)_{k\ge1}$ be independent and identically distributed positive
random variables, independent of the entire movement process, such that
\begin{equation}\label{eq:basic-moments}
 \E Z_1=1,\qquad \xi_k=Z_k-1,\qquad
 0<\sigma^2=\E\xi_1^2<\infty.
\end{equation}
Define the position with random step sizes and its centered noise by
\begin{equation}\label{eq:noise}
 H_n=\sum_{k=1}^n I_k Z_k,\qquad
 R_n=H_n-K_n=\sum_{k=1}^n I_k\xi_k.
\end{equation}
In particular, $\E(H_n\mid\mathcal F_n)=K_n$. The memory concerns
$(I_k)$, while the step sizes are sampled independently. We consider
\begin{equation}\label{eq:normalizations}
 V_n^2=\sum_{k=1}^n I_k\xi_k^2,\qquad
 U_n=\frac{R_n}{\sigma\sqrt{K_n}},\qquad
 T_n=\begin{cases}R_n/V_n,&V_n>0,\\0,&V_n=0.\end{cases}
\end{equation}
The statistic $U_n$ uses the conditional variance of $R_n$, whereas
$T_n$ is self-normalized. The convention at $V_n=0$ will be used
throughout. Neither $Z_1-1\eqdist1-Z_1$ nor any other symmetry
condition is assumed.

For $u>0$, put
\begin{equation}\label{eq:clock-moments}
 D_n(u)=\E K_n^{-u},\qquad
 h_n=\begin{cases}n^{-p/2},&\text{in case (A)},\\
                  n^{-1/2},&\text{in case (B)}.\end{cases}
\end{equation}

\subsection{Rates of normal approximation}

\quad We first present the normal approximation results for the
centered step-size noise. For real random variables $X,Y$, write
\[
 d_{\mathrm K}(X,Y)=\sup_{x\in\R}
       |\Pp(X\le x)-\Pp(Y\le x)|.
\]

\begin{theorem}\label{thm:normal}
Under \eqref{eq:basic-moments}, in either parameter set
\eqref{eq:regimes},
\begin{equation}\label{eq:clt}
 U_n\dto G,\qquad T_n\dto G.
\end{equation}
If, in addition, $\E|\xi_1|^{2+\delta}<\infty$ for some
$0<\delta\le1$, then, for $n\ge3$,
\begin{equation}\label{eq:be}
 d_{\mathrm K}(U_n,G)+d_{\mathrm K}(T_n,G)
 \le C D_n(\delta/2)\le C h_n^\delta.
\end{equation}
Moreover, $\E U_n=0$ and $\E U_n^2=1$ for every $n$.
\end{theorem}

We next consider the Wasserstein and Zolotarev distances.
For $s\ge1$, the Wasserstein distance is
\[
 W_s(X,Y)=\inf\bigl\{(\E|X'-Y'|^s)^{1/s}:
                X'\eqdist X,\ Y'\eqdist Y\bigr\}.
\]
For $0<s\le3$, let $\ell$ be the integer satisfying $\ell<s\le\ell+1$.
Let $\Lambda_s$ consist of all $\ell$-times continuously differentiable
functions $f:\R\to\R$ such that
\[
 |f^{(\ell)}(x)-f^{(\ell)}(y)|\le |x-y|^{s-\ell},
 \qquad x,y\in\R,
\]
where $f^{(0)}=f$. The Zolotarev distance is
\[
 \zeta_s(X,Y)=\sup_{f\in\Lambda_s}|\E f(X)-\E f(Y)|.
\]
It is finite under finite absolute $s$th moments and matching moments
of orders $1,\ldots,\ell$. The matching conditions cannot in general
be omitted, since $\Lambda_s$ contains polynomials of degree at most
$\ell$ with arbitrary coefficients.

\begin{theorem}\label{thm:distances}
Suppose that \eqref{eq:basic-moments} holds and $0<Z_1\le L$ almost
surely for some finite $L$. Put
\begin{equation}\label{eq:standardized-t}
 \mu_n=\E T_n,\qquad \tau_n^2=\Var(T_n),\qquad
 T_n^*=\frac{T_n-\mu_n}{\tau_n}.
\end{equation}
Then $\tau_n>0$ and, for $n\ge3$,
\begin{equation}\label{eq:t-moments}
 |\mu_n|+|\tau_n^2-1|\le C h_n.
\end{equation}
For every $0<s\le3$,
\begin{equation}\label{eq:zeta}
 \zeta_s(U_n,G)+\zeta_s(T_n^*,G)
 \le C_s h_n^{s\wedge1}.
\end{equation}
For every fixed $s\ge1$,
\begin{equation}\label{eq:w-general}
 W_s(U_n,G)+W_s(T_n,G)\le C_s D_n(s/2)^{1/s}.
\end{equation}
Consequently, in case~(A), the left-hand side of
\eqref{eq:w-general} is at most
\begin{equation}\label{eq:w-boundary}
 C_s\begin{cases}
 n^{-p/2},&1\le s<2,\\
 n^{-p/2}(\log n)^{1/2},&s=2,\\
 n^{-p/s},&s>2,
 \end{cases}
\end{equation}
and in case~(B) it is at most $C_s n^{-1/2}$.
\end{theorem}

\begin{remark}\label{rem:standardization}
The centering in \eqref{eq:standardized-t} is needed because
$\E T_n$ may be nonzero for a nonsymmetric step-size law, although
$\E R_n=0$.
Thus $\zeta_s(T_n,G)$ can be infinite for $s>1$; for $s>2$ the
variance must also match. No such adjustment is needed for $U_n$.
The high-order bounds in \eqref{eq:w-boundary} are upper bounds
obtained by mixing fixed-sample couplings; no optimality is asserted
for them.
\end{remark}

\begin{proposition}\label{prop:skewness}
Under the assumptions of Theorem~\ref{thm:distances}, if
$\E\xi_1^3\ne0$, then
\begin{equation}\label{eq:zeta-lower}
 \zeta_3(U_n,G)\ge
 \frac{|\E\xi_1^3|}{6\sigma^3}D_n(1/2)\asymp h_n.
\end{equation}
In particular, the order in \eqref{eq:zeta} for $U_n$ and $s=3$
is attained for such step-size laws.
\end{proposition}

\subsection{Self-normalized Cram\'er-type moderate deviations}

For $\varepsilon\in\{-1,1\}$ and $x\ge0$, set
\begin{equation}\label{eq:relative-error}
 \Delta_n^\varepsilon(x)=
 \left|\frac{\Pp(\varepsilon T_n\ge x)}{\Phibar(x)}-1\right|.
\end{equation}
All tail estimates are stated separately for the two signs, since
the centered step-size distribution need not be symmetric.

\begin{theorem}\label{thm:md-positive}
Assume case~(B), \eqref{eq:basic-moments}, and
$\E|\xi_1|^{2+\delta}<\infty$ for some $0<\delta\le1$.
There exist constants $c_0,C>0$ such that, for both choices of
$\varepsilon$ and all $n\ge3$,
\begin{equation}\label{eq:md-positive}
 \Delta_n^\varepsilon(x)
 \le C(1+x)^{2+\delta}n^{-\delta/2},
 \qquad 0\le x\le c_0 n^{\delta/[2(2+\delta)]}.
\end{equation}
In particular, for any nonnegative sequence
$x_n=o(n^{\delta/[2(2+\delta)]})$,
\begin{equation}\label{eq:md-positive-ratio}
 \sup_{0\le x\le x_n}\Delta_n^\varepsilon(x)\longrightarrow0.
\end{equation}
\end{theorem}

\begin{theorem}\label{thm:md-boundary}
Assume case~(A), \eqref{eq:basic-moments}, and
$\E|\xi_1|^{2+\delta}<\infty$ for some $0<\delta\le1$.
For both choices of $\varepsilon$, all $n\ge3$ and all $x\ge0$,
\begin{equation}\label{eq:md-boundary}
 \begin{split}
 \Delta_n^\varepsilon(x)\le C\bigl\{
 &(1+x)^{2+\delta}n^{-p\delta/2}\\
 &+(1+x)^{\,2(2+\delta)/\delta+1}
                  e^{x^2/2}n^{-p}\bigr\}.
 \end{split}
\end{equation}
Consequently, for every $0<\gamma<p$,
\begin{equation}\label{eq:log-window}
 \sup_{0\le x\le\sqrt{2\gamma\log n}}
       \Delta_n^\varepsilon(x)\longrightarrow0.
\end{equation}
The endpoint can be increased: for every $\eta>0$, the same uniform
convergence holds up to
\begin{equation}\label{eq:log-window-refined}
 x_n=\sqrt{2p\log n-
       \left(\frac{2(2+\delta)}\delta+1+\eta\right)\log\log n},
\end{equation}
for all sufficiently large $n$.
More precisely, if $x_n=o(\sqrt{\log n})$ is nonnegative, then, for
all sufficiently large $n$, uniformly for $0\le x\le x_n$,
\begin{equation}\label{eq:log-relative}
 \left|\log\frac{\Pp(\varepsilon T_n\ge x)}{\Phibar(x)}\right|
 \le C(1+x)^{2+\delta}n^{-p\delta/2}.
\end{equation}
\end{theorem}

Theorems~\ref{thm:md-positive} and \ref{thm:md-boundary} show
how small movement counts affect the range of normal tail
approximation. A finite third moment gives $x=o(n^{1/6})$ in
case~(B). The following example shows that a polynomial range
cannot hold for all bounded step-size laws in case~(A).

\begin{proposition}\label{prop:counterexample}
Let $q=0$, $p=1/2$, and let
\begin{equation}\label{eq:counter-law}
 \Pp(Z_1=9/8)=4/5,\qquad \Pp(Z_1=1/2)=1/5.
\end{equation}
Then $\E Z_1=1$ and $\sigma^2=1/16$. For
$m_n=\lfloor4\log n\rfloor$ and $x_n=\sqrt{m_n}$,
\begin{equation}\label{eq:counter-tail}
 \frac{\Pp(T_n\ge x_n)}{\Phibar(x_n)}\longrightarrow\infty.
\end{equation}
Thus no conclusion of the form
$\sup_{0\le x\le x_n}\Delta_n^+(x)\to0$ for every
$x_n=o(n^\eta)$, with a fixed $\eta>0$, can hold over all bounded
step-size laws in case~(A).
\end{proposition}

\begin{remark}
The logarithmic window in \eqref{eq:log-window} is a uniform sufficient
range under the stated assumptions. Proposition~\ref{prop:counterexample}
does not identify the optimal constant multiplying $\sqrt{\log n}$
for a given step law. The full-position results below use a different normalization,
because $H_n-\E H_n$ also contains $K_n-\E K_n$.
\end{remark}

\subsection{The full position}\label{subsec:full-main}

\quad We now turn to the full position. Put $a=p-q$,
$\rho=q/(1-a)$ and $m_n=\E H_n=\E K_n$. Then
\begin{equation}\label{eq:full-mean}
 m_n=\rho n+(1-\rho)\prod_{j=1}^{n-1}(1+a/j).
\end{equation}
The product is one at $n=1$ and also covers $a=-1$; for $a>-1$
it equals $\Gamma(n+a)/(\Gamma(1+a)\Gamma(n))$.
For $2\le k\le n$, define
\begin{equation}\label{eq:full-innovations}
 \theta_k=q+aK_{k-1}/(k-1),\qquad
 w_{n,k}=\prod_{j=k}^{n-1}(1+a/j).
\end{equation}
With $\mathcal G_k=\sigma(I_j,Z_j:j\le k)$ and trivial
$\mathcal G_0$, the array
\begin{equation}\label{eq:full-array}
 Y_{n,1}=\xi_1,\qquad
 Y_{n,k}=w_{n,k}(I_k-\theta_k)+I_k\xi_k\quad(k\ge2)
\end{equation}
consists of martingale differences and satisfies
\begin{equation}\label{eq:full-sum}
 H_n-m_n=\sum_{k=1}^nY_{n,k},\qquad
 s_n^2:=\Var(H_n)=\Var(K_n)+\sigma^2m_n>0.
\end{equation}
Define the standardized and self-normalized statistics by
\begin{equation}\label{eq:full-statistics}
 \mathcal U_n=(H_n-m_n)/s_n,\qquad
 \mathcal V_n^2=\sum_{k=1}^nY_{n,k}^2,\qquad
 \mathcal T_n=(H_n-m_n)/\mathcal V_n,
\end{equation}
with $\mathcal T_n=0$ when $\mathcal V_n=0$.
The denominator $\mathcal V_n$ accounts for both the movement
innovations and the centered step-size terms. It therefore differs
from $V_n$, which normalizes only $R_n$.

\begin{theorem}[Full-position limit laws]\label{thm:full-limits}
Assume \eqref{eq:basic-moments}.
In case~(B), $H_n/n\to\rho$ almost surely. For $p<1$ the following hold:
\begin{enumerate}
\item If $a<1/2$, then
\begin{equation}\label{eq:full-diffusive}
 s_n^2\sim n\left(\rho\sigma^2+\frac{\rho(1-\rho)}{1-2a}\right),\qquad
 \frac{H_n-m_n}{\sqrt n}\dto
 N\left(0,\rho\sigma^2+\frac{\rho(1-\rho)}{1-2a}\right).
\end{equation}
\item If $a=1/2$, then
\begin{equation}\label{eq:full-critical}
 s_n^2\sim\rho(1-\rho)n\log n,\qquad
 \frac{H_n-m_n}{\sqrt{n\log n}}\dto N(0,\rho(1-\rho)).
\end{equation}
In both these regimes,
$\mathcal U_n\dto G$, $\mathcal V_n^2/s_n^2\to1$ in probability,
and $\mathcal T_n\dto G$.
\item If $a>1/2$, there is a movement-measurable random variable $L_a$
such that
\begin{equation}\label{eq:full-super}
 \frac{H_n-\rho n}{n^a}\longrightarrow L_a
 \quad\hbox{almost surely and in }L^2,
 \qquad \E L_a=\frac{1-\rho}{\Gamma(1+a)}.
\end{equation}
Thus $(H_n-m_n)/n^a\to L_a-\E L_a$ in the same senses.
After centering by the random first-order limit, a Gaussian
second-order fluctuation remains:
\begin{equation}\label{eq:full-super-tail}
 \frac{H_n-\rho n-n^aL_a}{\sqrt n}\dto
 N\left(0,\rho\sigma^2+\frac{\rho(1-\rho)}{2a-1}\right).
\end{equation}
\end{enumerate}
If $p=1$ in case~(B), then $K_n=n$ almost surely and
$H_n=\sum_{k=1}^n Z_k$; consequently $s_n^2=n\sigma^2$ and
the ordinary i.i.d. conclusions apply for every $q>0$.

In case~(A), there is a positive Mittag--Leffler random variable $W_p$
such that
\begin{equation}\label{eq:full-boundary}
 \frac{H_n}{n^p}\longrightarrow W_p
 \quad\hbox{almost surely and in }L^2,\qquad
 \E W_p^r=\frac{r!}{\Gamma(1+pr)}\quad(r=1,2,\ldots).
\end{equation}
The moment formula gives $\Var(W_p)>0$ and
\begin{equation}\label{eq:full-boundary-standard}
 s_n^2\sim \Var(W_p)n^{2p},\qquad
 \mathcal U_n\dto
 \frac{W_p-1/\Gamma(1+p)}{\sqrt{\Var(W_p)}}.
\end{equation}
The limit in \eqref{eq:full-boundary-standard} is non-Gaussian.
Consequently, a Gaussian relative-tail approximation for
$\mathcal U_n$ cannot hold in general, even for fixed tail arguments.
\end{theorem}

The quantitative full-position results below concern $q>0$, $p<1$
and $a\le1/2$. Define, for $n\ge3$,
\begin{equation}\label{eq:full-rate-scales}
 e_n=\begin{cases}
 n^{-1/2+\max(a,0)},&a<1/2,\\
 (\log n)^{-1/2},&a=1/2,
 \end{cases}\qquad
 d_n=\begin{cases}
 n^{-1/2},&a<1/4,\\
 n^{-1/2}\sqrt{\log n},&a=1/4,\\
 n^{-(1-2a)},&1/4<a<1/2,\\
 (\log n)^{-1},&a=1/2.
 \end{cases}
\end{equation}
For $r>1$, set
\begin{equation}\label{eq:full-error-rate}
 \beta_n(r)=\left\{
 \frac{n+\sum_{k=2}^n w_{n,k}^{2r}}{s_n^{2r}}+d_n^r
 \right\}^{1/(2r+1)}.
\end{equation}

\begin{theorem}[Full-position normal approximation]\label{thm:full-be}
Suppose $q>0$, $p<1$ and $a\le1/2$.
If $\E|\xi_1|^{2r}<\infty$ for some $r>1$, then
\begin{equation}\label{eq:full-be}
 d_{\mathrm K}(\mathcal U_n,G)+d_{\mathrm K}(\mathcal T_n,G)
 \le C_r\beta_n(r).
\end{equation}
For the standardized position, the stronger nonuniform bound holds:
\begin{equation}\label{eq:full-nonuniform}
 |\Pp(\mathcal U_n\le x)-\Phi(x)|
 \le \frac{C_r\beta_n(r)}{1+|x|^{2r}},\qquad x\in\R.
\end{equation}
Under $\E|\xi_1|^{2+\delta}<\infty$, $0<\delta\le1$, the
right-hand side of \eqref{eq:full-be} may be replaced by
$C\min\{\beta_n(1+\delta/2),e_n^{\delta(2-\delta)/4}\}$.

Under the same $2r$th-moment assumption, for $1\le s<2r$,
\begin{equation}\label{eq:full-wasserstein}
 W_s(\mathcal U_n,G)\le C_{r,s}\beta_n(r)^{1/s}.
\end{equation}
For $0<s\le3$ with $s<2r$,
\begin{equation}\label{eq:full-zolotarev}
 \zeta_s(\mathcal U_n,G)\le C_{r,s}\beta_n(r)^{s\wedge1}.
\end{equation}
For bounded step sizes, $r>\max\{1,s/2\}$ can be chosen freely.
Here $\mathcal U_n$ has exactly mean zero and variance one.
These bounds hold for all sufficiently large $n$.
If $a=0$, independence improves \eqref{eq:full-be} to
$Cn^{-\delta/2}$ under a $(2+\delta)$th moment. For bounded step sizes,
$W_s(\mathcal U_n,G)\le C_sn^{-1/2}$ and
$\zeta_s(\mathcal U_n,G)\le C_sn^{-(s\wedge1)/2}$.
\end{theorem}

\begin{remark}
With only a third moment, \eqref{eq:full-be} gives order $n^{-1/8}$
for $a<1/3$, $n^{-1/8}(\log n)^{1/4}$ for $a=1/3$,
$n^{-3(1-2a)/8}$ for $1/3<a<1/2$, and $(\log n)^{-3/8}$ at
criticality. For bounded step sizes, choosing a sufficiently large fixed $r$
gives, for each $\eta>0$,
\begin{equation}\label{eq:full-bounded-rate}
 d_{\mathrm K}(\mathcal U_n,G)+d_{\mathrm K}(\mathcal T_n,G)
 \le C_\eta
 \begin{cases}
 n^{-\min\{1/4,\,1/2-a\}+\eta},&a<1/2,\\
 (\log n)^{-1/2+\eta},&a=1/2.
 \end{cases}
\end{equation}
The optimality of these rates for the present model is not asserted.
\end{remark}

\begin{theorem}[Full-position self-normalized moderate deviations]
\label{thm:full-md}
Suppose $q>0$, $p<1$, $a\le1/2$ and
$\E|\xi_1|^{2+\delta}<\infty$, $0<\delta\le1$.
For both $\varepsilon\in\{-1,1\}$,
\begin{equation}\label{eq:full-md-bound}
\begin{split}
 \left|\log\frac{\Pp(\varepsilon\mathcal T_n\ge x)}{\Phibar(x)}\right|
 \le C\bigg\{&x^{2+\delta}e_n^\delta+x^3d_n\\
 &+(1+x)\left(
 \frac{e_n^{\delta(2-\delta)/4}}{1+x^{\delta(2+\delta)/4}}
             +d_n(1+|\log d_n|)\right)\bigg\}.
\end{split}
\end{equation}
The estimate holds uniformly on $0\le x\le x_n$ for any
$x_n=o(e_n^{-1})$, for all sufficiently large $n$.
In particular,
\begin{equation}\label{eq:full-md-range}
 \sup_{0\le x\le x_n}
 \left|\frac{\Pp(\varepsilon\mathcal T_n\ge x)}{\Phibar(x)}-1\right|
 \longrightarrow0
\end{equation}
whenever
\begin{equation}\label{eq:full-md-explicit}
 x_n=\begin{cases}
 o\left(n^{(1/2-\max(a,0))\delta/(2+\delta)}\right),&a<1/2,\\
 o\left((\log n)^{\delta/[2(2+\delta)]}\right),&a=1/2.
 \end{cases}
\end{equation}
For a finite third moment, these ranges are respectively
$o(n^{(1/2-\max(a,0))/3})$ and $o((\log n)^{1/6})$.
For $p=1$, the i.i.d. range $o(n^{\delta/[2(2+\delta)]})$
holds instead, with $\mathcal T_n=T_n$.
\end{theorem}

\begin{remark}[Centering]\label{rem:full-scope}
The fluctuation of $K_n$ must be included when passing from $R_n$
to $H_n-m_n$. For $a>1/2$, the Gaussian limit in
\eqref{eq:full-super-tail} requires a random center. For $q=0$,
the limit in \eqref{eq:full-boundary-standard} is non-Gaussian.
These regimes are outside the scope of
Theorems~\ref{thm:full-be} and \ref{thm:full-md}.
\end{remark}

\section{Preliminary lemmas}
\label{sec:clock}

\subsection{The independent random-sum representation}

Let $(\widehat\xi_j)_{j\ge1}$ be independent copies of $\xi_1$,
independent of $(K_n)_{n\ge1}$, and define
\begin{equation}\label{eq:iid-sums}
 S_m=\sum_{j=1}^m\widehat\xi_j,\qquad
 Q_m=\sum_{j=1}^m\widehat\xi_j^2,\qquad
 A_m=\frac{S_m}{\sigma\sqrt m},\qquad
 B_m=\begin{cases}S_m/\sqrt{Q_m},&Q_m>0,\\0,&Q_m=0.\end{cases}
\end{equation}
Here $Q_m$ is the sum of squared centered step sizes in the independent
sample.

\begin{lemma}\label{lem:representation}
The pairs of processes $(R_n,V_n^2)_{n\ge1}$ and
$(S_{K_n},Q_{K_n})_{n\ge1}$ have the same finite-dimensional
distributions jointly with $(K_n)_{n\ge1}$. In particular, whenever
$\Pp(K_n=m)>0$,
\begin{equation}\label{eq:conditional-law}
 \Law\bigl((U_n,T_n)\mid K_n=m\bigr)=\Law(A_m,B_m).
\end{equation}
\end{lemma}

\begin{manuscriptproof}
Given all indicators, the centered step sizes at successive movement times are
i.i.d. with law $\Law(\xi_1)$. The sums at time $n$ contain their first
$K_n$ terms. This identifies all finite-dimensional conditional laws;
averaging proves the lemma.

\end{manuscriptproof}

\subsection{Movement-count estimates at the boundary}

\begin{lemma}\label{lem:clock-boundary}
In case~(A), $K_n\to\infty$ almost surely. Moreover,
\begin{align}
 \E K_n&=\frac{\Gamma(n+p)}{\Gamma(1+p)\Gamma(n)},
       \label{eq:clock-mean}\\
 \Pp(K_n=1)&=\frac{\Gamma(n-p)}{\Gamma(1-p)\Gamma(n)},
       \label{eq:clock-one}\\
 \sup_{1\le m\le n}\Pp(K_n=m)&\le C_p n^{-p}.
       \label{eq:clock-local}
\end{align}
For every $u>0$ and $n\ge3$,
\begin{equation}\label{eq:negative-boundary}
 D_n(u)\le C_{p,u}
 \begin{cases}
 n^{-pu},&0<u<1,\\
 n^{-p}\log n,&u=1,\\
 n^{-p},&u>1.
 \end{cases}
\end{equation}
In addition, $D_n(1/2)\asymp n^{-p/2}$.
\end{lemma}

\begin{manuscriptproof}
The transition probability is $pK_n/n$. Iterating the mean recursion
and the probability of repeated stops gives \eqref{eq:clock-mean}
and \eqref{eq:clock-one}. A bounded nondecreasing integer sequence
is eventually constant, whereas the conditional probability of staying
forever at any fixed $k$ after time $N$ is
$\prod_{j=N}^\infty(1-pk/j)=0$. Thus $K_n\to\infty$ almost surely.

For $0<u<1$, let $f_u(k)=\Gamma(k-u)/\Gamma(k)\asymp k^{-u}$.
Since $f_u(k+1)=(1-u/k)f_u(k)$, iteration yields
\begin{equation}\label{eq:gamma-moment}
 \E f_u(K_n)=\frac{\Gamma(1-u)\Gamma(n-pu)}
                  {\Gamma(1-pu)\Gamma(n)}\le Cn^{-pu}.
\end{equation}
For the remaining powers, put $h_{n,m}=\Pp(K_n=m)$, zero outside
$1\le m\le n$. The recursion is
\begin{equation}\label{eq:clock-recursion}
 h_{n+1,m}=(1-pm/n)h_{n,m}+p(m-1)h_{n,m-1}/n\quad(1\le m\le n),
\end{equation}
and $h_{n+1,n+1}=p h_{n,n}$. The two coefficients sum to $1-p/n$;
eventually $1-p/n\ge p$, so the maximum point mass is $O(n^{-p})$.
Therefore $D_n(u)\le Cn^{-p}\sum_{m=1}^n m^{-u}$ for $u\ge1$.
Finally $D_n(1/2)\ge(\E K_n)^{-1/2}\ge cn^{-p/2}$ by Jensen.
\end{manuscriptproof}

\subsection{Movement-count estimates for \texorpdfstring{$q>0$}{q>0}}

\begin{lemma}\label{lem:clock-positive}
In case~(B), there are constants $c_1,c_2,C>0$ such that
\begin{equation}\label{eq:clock-exponential}
 \Pp(K_n<c_1 n)\le C e^{-c_2 n}.
\end{equation}
Consequently, $K_n\to\infty$ almost surely and, for every $u>0$,
\begin{equation}\label{eq:negative-positive}
 D_n(u)\le C_u n^{-u},\qquad
 D_n(1/2)\asymp n^{-1/2}.
\end{equation}
\end{lemma}

\begin{manuscriptproof}
Realize $I_{k+1}=\ind\{W_{k+1}\le q+(p-q)K_k/k\}$ with
independent uniform $W_k$. On $E_n=\{K_n<qn/16\}$, for
$\lceil n/2\rceil\le k\le n-1$,
\[
 K_k/k<q/8,\qquad
 q+(p-q)K_k/k=q(1-K_k/k)+pK_k/k\ge q/2.
\]
Hence
\[
 E_n\subseteq\left\{
 \sum_{k=\lceil n/2\rceil}^{n-1}\ind\{W_{k+1}\le q/2\}<qn/16
 \right\}.
\]
A binomial Chernoff bound proves \eqref{eq:clock-exponential} with
$c_1=q/16$. Borel--Cantelli gives $K_n\to\infty$, and
$D_n(u)\le(c_1n)^{-u}+\Pp(E_n)\le C_un^{-u}$.
The lower bound for $D_n(1/2)$ follows from $K_n\le n$.
\end{manuscriptproof}

\section{Proofs of the normal approximation results}
\label{sec:normal-proofs}

\subsection{Fixed-sample estimates}

We first collect the estimates for independent sums used in the proofs.
Under \eqref{eq:basic-moments}, the ordinary central limit theorem
and the law of large numbers imply
\begin{equation}\label{eq:iid-clt}
 A_m\dto G,\qquad B_m\dto G.
\end{equation}
If $\E|\xi_1|^{2+\delta}<\infty$, $0<\delta\le1$, the
Lyapunov form of the Berry--Esseen inequality (see
Petrov \cite{Petrov1975}) and the self-normalized Berry--Esseen
bound of Shao \cite{Shao2005} give
\begin{equation}\label{eq:iid-be}
 d_{\mathrm K}(A_m,G)+d_{\mathrm K}(B_m,G)
 \le C_\delta\frac{\E|\xi_1|^{2+\delta}}{\sigma^{2+\delta}}
                     m^{-\delta/2}.
\end{equation}
The convention for $B_m$ when $Q_m=0$
can change a distribution function by at most
$\Pp(\xi_1=0)^m$; this exponentially small quantity is absorbed by
\eqref{eq:iid-be}.

\begin{lemma}\label{lem:iid-coupling}
Suppose $0<Z_1\le L$ almost surely and \eqref{eq:basic-moments}
holds. For every fixed $r\ge1$,
\begin{align}
 \|A_m\|_r+\|B_m\|_r&\le C_r,\label{eq:iid-moments}\\
 \|B_m-A_m\|_r&\le C_r m^{-1/2},\label{eq:iid-normalizer}\\
 W_r(A_m,G)+W_r(B_m,G)&\le C_r m^{-1/2}.
       \label{eq:iid-w}
\end{align}
For $0<s\le3$,
\begin{equation}\label{eq:iid-zeta}
 \zeta_s(A_m,G)\le C_s m^{-(s\wedge1)/2}.
\end{equation}
\end{lemma}

\begin{manuscriptproof}
Set $L_m=Q_m/(m\sigma^2)$. Hoeffding's inequality
\cite{Hoeffding1963} gives, for every fixed $r\ge1$,
\begin{equation}\label{eq:square-concentration}
 \|A_m\|_r\le C_r,\qquad \|L_m-1\|_r\le C_rm^{-1/2},\qquad
 \Pp(L_m<1/2)\le Ce^{-cm}.
\end{equation}
On $\{L_m\ge1/2\}$, $|B_m-A_m|\le C|A_m||L_m-1|$;
on its complement both variables are bounded by $C\sqrt m$.
H\"older's inequality proves \eqref{eq:iid-normalizer} and
\eqref{eq:iid-moments}. Rio \cite[Theorem 2.1]{Rio2009} gives
$W_r(A_m,G)\le C_rm^{-1/2}$ under bounded step sizes, so the triangle
inequality gives \eqref{eq:iid-w}.

For $s\le1$, $\zeta_s\le W_1^s$. For $1<s\le3$, subtract
from $f\in\Lambda_s$ its Taylor polynomial at zero of degree $\ell$;
its expectations match under $A_m,G$. The remainder satisfies
\begin{equation}\label{eq:test-growth}
 |f'(x)|\le C_s|x|^{s-1},\qquad
 |f(x)-f(y)|\le C_s|x-y|(|x|^{s-1}+|y|^{s-1}).
\end{equation}
Cauchy--Schwarz under an optimal $W_2$ coupling and the uniform
moment bounds yield \eqref{eq:iid-zeta}.
\end{manuscriptproof}

\subsection{Normal approximation}

\begin{manuscriptproof}[Proof of Theorem~\ref{thm:normal}]
By \eqref{eq:conditional-law},
$\E(U_n\mid K_n)=0$ and $\E(U_n^2\mid K_n)=1$.
Since $K_n\to\infty$ almost surely, bounded convergence transfers
\eqref{eq:iid-clt} to $U_n,T_n$. Averaging \eqref{eq:iid-be} gives
the bound $CD_n(\delta/2)$, and the clock lemmas give
\eqref{eq:be}.
\end{manuscriptproof}

\subsection{Moment distances}

\begin{manuscriptproof}[Proof of Theorem~\ref{thm:distances}]
Since $\xi_1$ is centered and nondegenerate, it takes both positive
and negative values with positive probability. The events on which
all selected centered step sizes have the same sign therefore show
that $\tau_n>0$. Mixing optimal fixed-sample couplings preserves
the Gaussian marginal and gives
\begin{equation}\label{eq:mixing-w}
 W_s(U_n,G)^s+W_s(T_n,G)^s\le C_sD_n(s/2).
\end{equation}
Thus \eqref{eq:w-general} and its explicit rates follow from the
clock estimates. Similarly, \eqref{eq:iid-normalizer} and
\eqref{eq:iid-moments} give
\begin{equation}\label{eq:mean-variance}
 |\mu_n|+|\tau_n^2-1|\le CD_n(1/2)\le Ch_n.
\end{equation}
Consequently $|\tau_n^{-1}-1|\le Ch_n$ for large $n$;
enlarge constants for the remaining values.

Averaging \eqref{eq:iid-zeta} yields
$\zeta_s(U_n,G)\le C_sD_n((s\wedge1)/2)\le C_sh_n^{s\wedge1}$.
Conditional on $K_n=m$, couple $T_n^*$ and $U_n$ as
$(\tau_n^{-1}(B_m-\mu_n),A_m)$. For each fixed $r\ge1$,
\begin{equation}\label{eq:affine-coupling}
 \|\tau_n^{-1}(B_m-\mu_n)-A_m\|_r\le C_r(m^{-1/2}+h_n),
\end{equation}
and both variables have uniformly bounded $r$th moments.
For $s\le1$, this bounds $\zeta_s(T_n^*,U_n)$ by
$C_s\{D_n(s/2)+h_n^s\}$. For $1<s\le3$, subtract the Taylor
polynomial in the unconditional expectation, where moments match.
Using \eqref{eq:test-growth} and Cauchy--Schwarz conditionally on $K_n=m$,
then averaging \eqref{eq:affine-coupling}, gives $C_sh_n$.
The triangle inequality proves \eqref{eq:zeta}.
\end{manuscriptproof}

\begin{manuscriptproof}[Proof of Proposition~\ref{prop:skewness}]
Independence and centering of the step-size noise give
$\E A_m^3=(\E\xi_1^3)/(\sigma^3\sqrt m)$. Consequently,
\[
 \E U_n^3=\frac{\E\xi_1^3}{\sigma^3}D_n(1/2).
\]
The function $f(x)=x^3/6$ belongs to $\Lambda_3$, because
$f''(x)=x$ has Lipschitz constant one. Since $\E G^3=0$,
the definition of $\zeta_3$ gives \eqref{eq:zeta-lower}.
The clock lemmas and \eqref{eq:zeta} complete the proof.
\end{manuscriptproof}

\section{Proofs of the moderate-deviation results}
\label{sec:md-proofs}

\subsection{The fixed-sample input and its transfer}

Under $\E|\xi_1|^{2+\delta}<\infty$, $0<\delta\le1$, the
finite-moment relative tail estimate of Jing, Shao and Wang
\cite{Jing2003}, applied to $\varepsilon\widehat\xi_j$, gives
constants $c_*,C>0$ such that
\begin{equation}\label{eq:iid-md}
 \left|\frac{\Pp(\varepsilon B_m\ge x)}{\Phibar(x)}-1\right|
 \le C(1+x)^{2+\delta}m^{-\delta/2},
 \quad 0\le x\le c_*m^{\delta/[2(2+\delta)]}.
\end{equation}
The event $Q_m=0$ has probability $\Pp(\xi_1=0)^m$, so our convention
does not affect this bound. Conditioning on clock values
$\mathcal A$ where \eqref{eq:iid-md} applies gives
\begin{equation}\label{eq:split-tail}
 \Delta_n^\varepsilon(x)\le
 C(1+x)^{2+\delta}\E\{K_n^{-\delta/2};K_n\in\mathcal A\}
 +\frac{2\Pp(K_n\notin\mathcal A)}{\Phibar(x)}.
\end{equation}
We use $\Phibar(x)^{-1}\le C(1+x)e^{x^2/2}$ for $x\ge0$.

\subsection{The case \texorpdfstring{$q>0$}{q>0}}

\begin{manuscriptproof}[Proof of Theorem~\ref{thm:md-positive}]
Take $\mathcal A=\{m\ge c_1n\}$ in \eqref{eq:split-tail}.
For sufficiently small $c_0$, the fixed-sample estimate applies
throughout $x\le c_0n^{\delta/[2(2+\delta)]}$, and
\[
 \Delta_n^\varepsilon(x)\le
 C(1+x)^{2+\delta}n^{-\delta/2}
 +C(1+x)e^{-c_2n+x^2/2}.
\]
Since $x^2=o(n)$ uniformly in this range, the second term is
absorbed by the first. The uniform convergence follows immediately.
\end{manuscriptproof}

\subsection{The case \texorpdfstring{$q=0$}{q=0}}

\begin{manuscriptproof}[Proof of Theorem~\ref{thm:md-boundary}]
Let $M_x=\lceil A(1+x)^{2(2+\delta)/\delta}\rceil$, with $A$
large enough that \eqref{eq:iid-md} applies whenever $m\ge M_x$.
The local clock bound implies
\begin{equation}\label{eq:small-clock}
 \Pp(K_n<M_x)\le Cn^{-p}(1+x)^{2(2+\delta)/\delta}.
\end{equation}
Use \eqref{eq:split-tail} with $\mathcal A=\{m\ge M_x\}$,
$D_n(\delta/2)\le Cn^{-p\delta/2}$ and the normal-tail bound.
This proves \eqref{eq:md-boundary}, also when $M_x>n$.

On $x\le\sqrt{2\gamma\log n}$, $\gamma<p$, the two error terms
are polynomial logarithmic factors times $n^{-p\delta/2}$ and
$n^{-p+\gamma}$, respectively. At the larger endpoint
\eqref{eq:log-window-refined}, the second is
$O((\log n)^{-\eta/2})$, and the first still tends to zero.
Finally, for $x_n=o(\sqrt{\log n})$, the second error term divided
by the first is at most $n^{-p(1-\delta/2)+o(1)}$.
The first tends uniformly to zero, so
$|\log(1+y)|\le2|y|$ for $|y|\le1/2$ proves
\eqref{eq:log-relative}.
\end{manuscriptproof}

\subsection{A bounded nonsymmetric counterexample}

\begin{manuscriptproof}[Proof of Proposition~\ref{prop:counterexample}]
Here $\xi_1=1/8$ with probability $4/5$ and $-1/2$ otherwise,
so $\E\xi_1=0$ and $\sigma^2=1/16$. When $p=1/2$, induction
in \eqref{eq:clock-recursion}, including its upper endpoint, gives
\begin{equation}\label{eq:half-clock-law}
 \Pp(K_n=m)=2^{m+1-2n}\binom{2n-m-1}{n-1},\qquad1\le m\le n.
\end{equation}
Stirling's formula yields
\begin{equation}\label{eq:half-clock-local}
 \Pp(K_n=m)=\frac{1+o(1)}{\sqrt{\pi n}},
 \qquad1\le m\le C\log n,
\end{equation}
uniformly for each fixed $C$: the binomial index is within
$O(\log n)$ of its center, with relative correction $1+o(1)$.
If all selected centered step sizes equal $1/8$, then
$T_n=\sqrt{K_n}$. Hence
\[
 \Pp(T_n\ge\sqrt m)\ge\Pp(K_n=m)(4/5)^m.
\]
For $m=m_n=\lfloor4\log n\rfloor$, the normal-tail asymptotic gives
\begin{equation}\label{eq:counter-divergence}
 \frac{\Pp(T_n\ge\sqrt{m_n})}{\Phibar(\sqrt{m_n})}
 \ge n^{3/2+4\log(4/5)+o(1)}\longrightarrow\infty.
\end{equation}
The exponent is positive and $\sqrt{m_n}=o(n^\eta)$ for every
$\eta>0$, proving both assertions.
\end{manuscriptproof}

\section{Proofs of the full-position results}\label{sec:full-proofs}

\subsection{Martingale representation and variance concentration}

Taking expectations in \eqref{eq:transition} gives
$m_{n+1}=(1+a/n)m_n+q$, proving \eqref{eq:full-mean}.
The centered recursion is
\[
 K_{n+1}-m_{n+1}=(1+a/n)(K_n-m_n)+I_{n+1}-\theta_{n+1}.
\]
Iteration proves \eqref{eq:full-sum}. The centered step-size terms are conditionally
centered and orthogonal to the movement innovations. Consequently
the predictable quadratic variation of the array is
\begin{equation}\label{eq:full-bracket}
 \mathcal B_n:=\sum_{k=1}^n\E(Y_{n,k}^2\mid\mathcal G_{k-1})
 =\sigma^2+\sum_{k=2}^n
   \{w_{n,k}^2\theta_k(1-\theta_k)+\sigma^2\theta_k\},
 \qquad \E\mathcal B_n=s_n^2.
\end{equation}
The product formula yields $\log w_{n,k}=a\log(n/k)+O(1/k)$,
uniformly for $n\ge k\ge2$. Thus
\begin{equation}\label{eq:full-weight-estimates}
 w_{n,k}\asymp(n/k)^a,\qquad
 \sum_{k=2}^n w_{n,k}^2\sim
 \begin{cases}n/(1-2a),&a<1/2,\\n\log n,&a=1/2.
 \end{cases}
\end{equation}

\begin{lemma}\label{lem:full-concentration}
For $q>0$, $p<1$ and $a\le1/2$, the variance asymptotics in
\eqref{eq:full-diffusive}--\eqref{eq:full-critical} hold, and
\begin{equation}\label{eq:full-concentration}
 \Pp(|\mathcal B_n/s_n^2-1|>t)\le Ce^{-ct^2/d_n^2},\qquad t>0.
\end{equation}
If $\E|\xi_1|^{2+\delta}<\infty$, $0<\delta\le1$, then
\begin{equation}\label{eq:full-conditional-moment}
 \E(|Y_{n,k}|^{2+\delta}\mid\mathcal G_{k-1})
 \le(Ce_ns_n)^\delta\E(Y_{n,k}^2\mid\mathcal G_{k-1}).
\end{equation}
Moreover $\Pp(\mathcal V_n=0)\le Ce^{-cn}$.
\end{lemma}

\begin{manuscriptproof}
Let $\|Z\|_{\psi_2}=\inf\{u>0:\E e^{Z^2/u^2}\le2\}$.
The conditional Hoeffding bound for
$K_j-m_j=\sum_{k=2}^j w_{j,k}(I_k-\theta_k)$, together with
\eqref{eq:full-mean} and \eqref{eq:full-weight-estimates}, gives
\begin{equation}\label{eq:full-clock-psi}
 \|K_j-\rho j\|_{\psi_2}\le
 \begin{cases}C\sqrt j,&a<1/2,\\
 C\sqrt{j\log(j+1)},&a=1/2.
 \end{cases}
\end{equation}
Replacing $\theta_k$ by $\rho$ in \eqref{eq:full-bracket}
has expectation error at most
$C\sum_{k=2}^n(w_{n,k}^2+1)/\sqrt{k-1}=o(n)$ for $a<1/2$,
and $C\sum_{k=2}^n(n/k+1)\sqrt{\log k/(k-1)}=O(n)$ at
criticality. This proves the variance asymptotics.

For concentration, reveal $I_2,\ldots,I_{n-1}$ successively.
The count transition kernels preserve stochastic order: from adjacent
counts $z,z+1$, the next counts have supports $\{z,z+1\}$ and
$\{z+1,z+2\}$. Under their monotone coupling, the expected count
gap at time $k-1$, starting with gap one at time $j$, is
$w_{k-1,j}$ (with $w_{j,j}=1$). This follows by iterating the
conditional mean recursion, also when $a<0$.
The $k$th summand of \eqref{eq:full-bracket}, as a function of
$K_{k-1}$, is Lipschitz with constant
$|a|(w_{n,k}^2+\sigma^2)/(k-1)$.
Hence the $j$th difference of the Doob martingale of $\mathcal B_n$
has conditional range length at most
\[
 c_{n,j}:=|a|\sum_{k=j+1}^n
 \frac{w_{n,k}^2+\sigma^2}{k-1}w_{k-1,j}
 \le C\begin{cases}1,&a<0,\\(n/j)^{2a},&a>0.\end{cases}
\]
For $a=0$, $\mathcal B_n$ is deterministic. For $a\ne0$,
conditional Hoeffding and \eqref{eq:full-weight-estimates} give
\[
 \|\mathcal B_n-\E\mathcal B_n\|_{\psi_2}^2
 \le C\sum_{j=2}^{n-1}c_{n,j}^2
 \le C\begin{cases}
 n,&a<1/4,\\ n\log n,&a=1/4,\\ n^{4a},&1/4<a\le1/2.
 \end{cases}
\]
Division by $s_n^4$ proves \eqref{eq:full-concentration}, including
the square-root logarithm at $a=1/4$.

For $k\ge2$, $|I_k-\theta_k|\le1$ and independence of the step sizes give
\begin{align*}
 \E(|Y_{n,k}|^{2+\delta}\mid\mathcal G_{k-1})
 &\le C\{w_{n,k}^{2+\delta}\theta_k(1-\theta_k)
                      +\theta_k\E|\xi_1|^{2+\delta}\}\\
 &\le C\max(w_{n,k},1)^\delta
       \E(Y_{n,k}^2\mid\mathcal G_{k-1}).
\end{align*}
Now $\max(1,\max_{k\ge2}w_{n,k})/s_n\le Ce_n$; the $k=1$
case is immediate. Finally let $r=\sup_z\Pp(\xi_1=z)<1$.
Given the indicators, $\mathcal V_n=0$ prescribes one step-size value
at each movement. Thus
$\Pp(\mathcal V_n=0)\le\E r^{K_n}\le Ce^{-cn}$ by
Lemma~\ref{lem:clock-positive}.
\end{manuscriptproof}

\subsection{Limit laws}

\begin{manuscriptproof}[Proof of Theorem~\ref{thm:full-limits}]
The movement strong law gives $K_n/n\to\rho$ for $q>0$
\cite{Coletti2019,Bercu2022}. At $a=-1$ the same conclusion follows
from \eqref{eq:full-clock-psi} and Borel--Cantelli.
The step sizes at successive movement times are i.i.d.; their strong law
gives $H_n/K_n\to1$, hence $H_n/n\to\rho$.

For $a\le1/2$ and $p<1$, Lemma~\ref{lem:full-concentration} gives
$\mathcal B_n/s_n^2\to1$ in probability, and
$\max_{k\ge2}w_{n,k}/s_n\to0$. The bounded movement terms
therefore satisfy the Lindeberg condition. For the centered step-size terms, it follows from
\[
 \frac1{s_n^2}\sum_{k=1}^n
 \E\{I_k\xi_k^2;|\xi_k|>\eta s_n\}
 \le C\E\{\xi_1^2;|\xi_1|>\eta s_n\}\longrightarrow0.
\]
Splitting $Y_{n,k}$ into its two terms verifies the condition for the
array; the additional movement-square term on a large step-size event is bounded
by $C\Pp(|\xi_1|>\eta s_n/2)$.
The martingale CLT proves $\mathcal U_n\dto G$.
The same Lindeberg condition, by truncation of the squared increments,
gives $(\mathcal V_n^2-\mathcal B_n)/s_n^2\to0$ in probability.
Slutsky's theorem proves the remaining assertions in these regimes.

For $a>1/2$, the known clock limit is
$(K_n-\rho n)/n^a\to L_a$ almost surely and in $L^2$
\cite{Bercu2022}. Since $\E R_n^2=O(n)$, $R_n$ vanishes at this
scale in $L^2$. The martingale series $\sum k^{-a}I_k\xi_k$ has
summable increment variances, so Kronecker's lemma gives the same
almost-sure conclusion. Formula \eqref{eq:full-mean} identifies $\E L_a$.
For the second-order limit, the clock theorem in \cite{Bercu2022} gives
\[
 (K_n-\rho n-n^aL_a)/\sqrt n\dto
 N(0,\rho(1-\rho)/(2a-1)).
\]
Writing $\varphi(t)=\E e^{it\xi_1}$, conditionally on all indicators,
\begin{equation}\label{eq:full-asymptotic-independence}
 \E(e^{itR_n/\sqrt n}\mid\mathcal F_\infty)
 =\varphi(t/\sqrt n)^{K_n}\longrightarrow e^{-\rho\sigma^2t^2/2}
 \quad\hbox{almost surely and in }L^1.
\end{equation}
Multiplication by the bounded characteristic-function factor of the
clock fluctuation proves asymptotic independence, so the variances add.
For $p=1$, induction gives $I_k=1$ and all claims reduce to i.i.d. sums.

For $q=0$, Miyazaki and Takei \cite{Miyazaki2020} give
$K_n/n^p\to W_p$ almost surely and in every fixed $L^r$, with the
moments in \eqref{eq:full-boundary}. Now $H_n/K_n\to1$ and
$\E(R_n/n^p)^2=\sigma^2m_n/n^{2p}=O(n^{-p})$.
This proves the almost-sure and $L^2$ limits for $H_n$ and, using
\eqref{eq:full-sum}, the variance asymptotic. The centered limit is
nondegenerate and bounded below, so it is not Gaussian.
\end{manuscriptproof}

\subsection{Moderate deviations}

\begin{manuscriptproof}[Proof of Theorem~\ref{thm:full-md}]
Apply Fan and Shao \cite[Theorem 2.1, Remark 2.1]{FanShao2025}
to $\pm Y_{n,k}$ with deterministic variance scale $s_n$.
Lemma~\ref{lem:full-concentration} verifies their Gaussian-concentration
condition (A1$'$) and conditional-moment condition (A2), with constant
multiples of $d_n,e_n$. Their bound is \eqref{eq:full-md-bound};
the $x^3d_n$ term corresponds to (A1$'$). Since $d_n\le Ce_n$,
the admissible range is $x=o(e_n^{-1})$. The zero-denominator
probability is exponentially small and negligible relative to the
normal tail there.

All terms vanish for
$x_n=o(\min\{e_n^{-\delta/(2+\delta)},d_n^{-1/3}\})$.
For the fractional term in \eqref{eq:full-md-bound}, we have
\[
 \frac{\delta(2-\delta)}4
 -\frac{\delta}{2+\delta}\left(1-\frac{\delta(2+\delta)}4\right)
 =\frac{\delta^2}{2(2+\delta)}>0.
\]
Equation~\eqref{eq:full-rate-scales} gives
$e_n^{-\delta/(2+\delta)}=O(d_n^{-1/3})$, proving
\eqref{eq:full-md-explicit}. For $p=1$, use \eqref{eq:iid-md}.
\end{manuscriptproof}

\subsection{Normal approximation and moment distances}

\begin{manuscriptproof}[Proof of Theorem~\ref{thm:full-be}]
For $\E|\xi_1|^{2r}<\infty$, $r>1$, the normalized array satisfies
\begin{equation}\label{eq:full-aggregate-moments}
 \sum_{k=1}^n\E|Y_{n,k}/s_n|^{2r}
 +\E|\mathcal B_n/s_n^2-1|^r
 \le C_r\left\{\frac{n+\sum_{k=2}^n w_{n,k}^{2r}}{s_n^{2r}}
                     +d_n^r\right\}.
\end{equation}
The first term follows directly from \eqref{eq:full-array};
the second follows by integrating \eqref{eq:full-concentration}.
The nonuniform martingale bound of Haeusler--Joos
\cite{HaeuslerJoos1988} applied to \eqref{eq:full-aggregate-moments}
proves \eqref{eq:full-nonuniform}. Together with Fan--Shao
\cite[Theorem 2.1]{FanShao2018} for the self-normalized sum,
it proves \eqref{eq:full-be}.
The zero-denominator convention adds only $Ce^{-cn}$.
Summing $(n/k)^{2r}$ gives the stated third-moment rates.
For bounded step sizes, as the fixed $r$ increases, the power of $n$ in
$\beta_n(r)$ approaches $-\min\{1/4,1/2-a\}$, and the critical
power of $\log n$ is $-r/(2r+1)$. Any logarithmic factors in the
diffusive case are absorbed by $n^\eta$, proving
\eqref{eq:full-bounded-rate}.

It remains to prove the alternative bound under a
$(2+\delta)$th moment. Put $\kappa=\delta(2-\delta)/4$.
Equation~\eqref{eq:full-md-bound}, multiplied by $\Phibar(x)$
for $x\le\sqrt{4\kappa|\log e_n|}$ and combined with monotonicity
beyond that point, gives $d_{\mathrm K}(\mathcal T_n,G)\le Ce_n^\kappa$;
here $d_n(1+|\log d_n|)=o(e_n^\kappa)$.
The martingale moment inequality for powers in $(1,2]$ and
\eqref{eq:full-conditional-moment} give
\[
 \E| (\mathcal V_n^2-\mathcal B_n)/s_n^2|^{1+\delta/2}
 \le Ce_n^\delta.
\]
Markov's inequality, \eqref{eq:full-concentration}, and
$\kappa(2+\delta/2)\le\delta$ show that
$\Pp(|\mathcal V_n/s_n-1|>Ce_n^\kappa)\le Ce_n^\kappa$.
Comparing distribution functions on the complementary event transfers
the same bound to $\mathcal U_n$.

Let $F_n$ be the distribution function of $\mathcal U_n$.
For $0<s<2r$, direct integration of \eqref{eq:full-nonuniform} gives
\begin{equation}\label{eq:full-cdf-integral}
 \int_{\R}|x|^{s-1}|F_n(x)-\Phi(x)|\,dx
 \le C_r\beta_n(r)\int_{\R}\frac{|x|^{s-1}}{1+|x|^{2r}}\,dx
 \le C_{r,s}\beta_n(r).
\end{equation}
For $s\ge1$, the monotone coupling inequality
$W_s^s\le s2^{s-1}\int|x|^{s-1}|F_n-\Phi|$ proves
\eqref{eq:full-wasserstein}. For $s\le1$, use $\zeta_s\le W_1^s$.
For $1<s\le3$, subtract the Taylor polynomial from $f\in\Lambda_s$;
its expectations match and $|f'(x)|\le C_s|x|^{s-1}$.
Integration by parts and \eqref{eq:full-cdf-integral} prove
\eqref{eq:full-zolotarev} directly.

If $a=0$, the array terms are $Y_{n,1}=\xi_1$ and
$Y_{n,k}=I_kZ_k-p$ for $k\ge2$, all independent, with the latter
terms identically distributed. The independent-sum Berry--Esseen
bounds give $Cn^{-\delta/2}$. With bounded step sizes, Rio's bound
\cite{Rio2009} applies to the last $n-1$ terms; adding the first
term and matching the exact variance costs $O(n^{-1/2})$ in each
fixed $W_s$. The Taylor argument in \eqref{eq:test-growth} gives
the stated Zolotarev rates.
\end{manuscriptproof}

\section{Concluding remarks}\label{sec:conclusion}

The random-sum representation yields quantitative estimates for
the centered step-size noise, while the weighted martingale array
includes the movement fluctuations in the full position. The
boundary case $q=0$ leads to a logarithmic moderate-deviation
range for the noise and a non-Gaussian limit for the full position.
Further questions concern optimal normal approximation rates,
the endpoint of the logarithmic range for a given step-size law,
and moderate deviations about random centers in the
superdiffusive regime.

\end{document}